\documentclass[11pt, reqno]{amsart}
\usepackage[utf8]{inputenc}
\usepackage{amsthm}
\usepackage{amsmath,amsfonts,amssymb}
\usepackage[a4paper, margin=1in]{geometry}
\usepackage{bm}
\usepackage{pgfplots}
\pgfplotsset{compat=1.18}
\usepgfplotslibrary{groupplots}

\usepackage{siunitx}
\usepackage{cite}
\usepackage{mwe}
\usepackage{graphicx}
\usepackage{caption}
\usepackage{subcaption}
\usepackage{hyperref}
\usepackage{soul}

\begin{document}

\title[Schr\"odinger Regularization of SWE]{%
A Schr\"odinger-Based Dispersive Regularization Approach for Numerical Simulation of One-Dimensional Shallow Water Equations
}

\author{Guosheng Fu}
\address{Department of Applied and Computational Mathematics and
Statistics, University of Notre Dame, USA.}
\email{gfu@nd.edu}
\author{Chun Liu}
\address{Department of Applied Mathematics, Illinois Institute of Technology, USA.}
\email{cliu124@iit.edu}
 \thanks{
G. Fu’s research is partially supported by NSF grant DMS-2410741.
Part of this research was performed while G. Fu was visiting the Institute for Mathematical and Statistical Innovation (IMSI), which is supported by the National Science Foundation (Grant No. DMS-2425650).
C. Liu’s research is partially supported by NSF grant DMS-2216926 and DMS-2410742.
}


% =========================
\begin{abstract}
We propose a novel dispersive regularization framework for the numerical simulation
of the one-dimensional shallow water equations (SWE).
The classical hyperbolic system is regularized by a third-order dispersive term
in the momentum equation, which renders the system equivalent, via the Madelung
transform, to a defocusing cubic nonlinear Schr\"odinger equation with a drift
term induced by bottom topography.

Instead of solving the shallow water equations directly, we solve the associated
Schr\"odinger equation and recover the hydrodynamic variables through a simple
postprocessing procedure.
This approach transforms the original nonlinear hyperbolic system into a
semilinear complex-valued equation, which can be efficiently approximated using
a Strang time-splitting method combined with a spectral element discretization
in space.

Numerical experiments demonstrate that, in subcritical regimes without shock
formation, the Schr\"odinger regularization provides an $O(\varepsilon)$ approximation
to the classical shallow water solution, where $\varepsilon$ denotes the
regularization parameter.
Importantly, we observe that this convergence behavior persists even in the
presence of moving wetting--drying interfaces, where vacuum states emerge and
standard shallow water solvers often encounter difficulties.
These results suggest that the Schr\"odinger-based formulation offers a robust
and promising alternative framework for the numerical simulation of shallow
water flows with dry states.
\end{abstract}

\maketitle

% =========================
\section{Introduction and Motivation}
% =========================

The shallow water equations (SWE) are a fundamental model for free-surface flows
arising in river hydraulics, coastal engineering, and geophysical fluid dynamics
\cite{toro2024computational}.
Despite their apparent simplicity, the numerical simulation of SWE remains
challenging, particularly in the presence of wetting--drying interfaces and
vacuum states, where the water depth vanishes and the equations degenerate.

A wide range of numerical methods has been developed for SWE, including finite
volume, finite difference, and discontinuous Galerkin (DG) schemes.
To handle dry states, these methods typically rely on carefully designed
positivity-preserving limiters, hydrostatic reconstructions, and ad hoc
regularizations of the velocity field
\cite{LeVeque02,Audusse04,Bunya09,XingShu11,Xing13}.
While such approaches are effective in many practical applications, they often
entail substantial algorithmic complexity and problem-dependent tuning, and the
robust treatment of moving wetting--drying fronts remains a delicate issue.

An alternative modeling approach is the use of \emph{dispersive regularizations}
of hyperbolic systems.
Rather than introducing artificial viscosity, dispersive terms are incorporated
in a way that preserves conservation laws and, in many cases, an underlying
Hamiltonian or variational structure.
Such dispersive mechanisms arise naturally in the hydrodynamic formulation of
the nonlinear Schr\"odinger equation and in Euler--Korteweg-type models
\cite{Madelung1927,CarlesDanchinSaut2012}.
Through the Madelung transform, solutions of the Schr\"odinger equation can be
interpreted in hydrodynamic variables, leading to a dispersively regularized
Euler-type system involving a quantum (or Bohm) potential.
These equations have been extensively studied in the context of the
semiclassical limit of Schr\"odinger dynamics
\cite{Grenier1998,Carles2008,BaoJinMarkowich2002}, with the primary objective of
understanding the behavior and limiting properties of the Schr\"odinger equation
itself, rather than serving as numerical regularizations of classical hyperbolic
conservation laws.

In this work, we build on these ideas and exploit the dispersive structure
induced by the Madelung transform as a numerical approximation framework for the
one-dimensional shallow water equations.
Specifically, we introduce a dispersive regularization of SWE that is equivalent
to a defocusing cubic nonlinear Schr\"odinger equation with a drift term induced
by bottom topography.
Rather than solving the shallow water equations directly, we solve the
associated Schr\"odinger equation and recover the physical variables through a
simple postprocessing procedure.
This reformulation converts the original nonlinear hyperbolic system into a
semilinear complex-valued equation, which is well suited for efficient numerical
approximation using time-splitting methods and high-order spatial
discretizations.

From a modeling and numerical perspective, the Schr\"odinger formulation has
several appealing features.
First, the water height is represented as the squared modulus of the wave
function, which guarantees non-negativity and naturally accommodates vacuum
states.
Second, the resulting semilinear structure allows for efficient numerical
schemes based on Strang splitting, in which the dispersive and potential
subproblems can be treated separately.
Finally, the dispersive regularization provides a controlled small-parameter
approximation to the dispersionless shallow water equations in regimes where the
solution remains smooth.

The main contributions of this paper can be summarized as follows:
\begin{itemize}
  \item We introduce a dispersive regularization of the one-dimensional shallow
        water equations that is equivalent to a defocusing cubic nonlinear
        Schr\"odinger equation with drift induced by bottom topography.
  \item We present an energetic variational (EnVarA) interpretation of the
        regularized formulation, clarifying the origin of the dispersive term
        from an energy-based viewpoint.
  \item We develop an efficient numerical method based on spectral element
        discretization in space and Strang splitting in time, leading to linear
        constant-coefficient subproblems in the dispersive step and analytic,
        pointwise updates in the potential step.
  \item Through a series of numerical experiments, we demonstrate that the
        proposed approach yields $\mathcal{O}(\varepsilon)$ accuracy in
        subcritical and shock-free regimes, including challenging scenarios
        involving moving wetting--drying interfaces and vacuum formation.
\end{itemize}

The remainder of the paper is organized as follows.
In Section~2, we introduce the dispersive regularization and its equivalence to a
nonlinear Schr\"odinger equation, and describe the numerical discretization.
Section~3 presents numerical results illustrating the strengths and limitations
of the proposed approach.
Concluding remarks and directions for future work are given in Section~4.

% =========================
\section{Schr\"odinger Regularization and Numerical Method}
% =========================

\subsection{Shallow Water Equations}

We consider the one-dimensional shallow water equations with bottom topography
$b(x)$,
\begin{subequations}    
\label{eq:SWE}
\begin{align}
\label{eq:SWE1}
h_t + (hu)_x &= 0, \\
\label{eq:SWE2}
(hu)_t + \left( hu^2 + \frac{1}{2} g h^2 \right)_x &= - g h b_x ,
\end{align}
\end{subequations}
where $g>0$ denotes the gravitational constant, $h(x,t)\ge0$ is the water depth,
and $u(x,t)$ is the depth-averaged velocity.

\subsection{Energetic variational approach and dispersive regularization}

The shallow water equations can be interpreted within the framework of the
energetic variational approach (EnVarA), which provides a systematic way to
derive governing equations from energy functionals by combining the least action
principle with force balance laws; see, e.g.,
\cite{Liu2009,Giga18}.
This perspective is particularly useful for elucidating the conservative
structure of the equations and for incorporating additional physical effects
through modifications of the energy.

For the one-dimensional shallow water equations, the total energy is given by
\[
\mathcal{E}(h,u)
=
\mathcal{K}(h,u)+\mathcal{F}(h)
=
\int_{\Omega} \frac{1}{2}h|u|^2\,dx
+
\int_{\Omega}\left(\frac{1}{2}g h^2 + g b(x)h\right)\,dx,
\]
where $\mathcal{K}$ denotes the kinetic energy, $\mathcal{F}$ the potential
energy associated with gravity and bottom topography, and $\Omega$ is the spatial domain.
The kinematic constraint relating the water height and the velocity is the mass
conservation law \eqref{eq:SWE1}.
In the absence of dissipation, smooth solutions satisfy the energy conservation
law
\begin{equation}
\label{ener0}
\frac{d}{dt}\mathcal{E}(h,u)=0.
\end{equation}

Within EnVarA, the conservative dynamics are derived by introducing a flow map
$x=x(X,t)$ satisfying
\[
\partial_t x(X,t)=u(x(X,t),t),\qquad x(X,0)=X,
\]
together with mass conservation. In one space dimension, mass conservation is
equivalently expressed in Lagrangian form as
\[
h(x(X,t),t)\,x_X(X,t)=h_0(X),
\]
or, in Eulerian form, as \eqref{eq:SWE1}.

Consider the kinetic-energy action functional in Lagrangian coordinates,
\[
\mathcal{A}[x]
=
\int_0^T \int_{\Omega_0} \frac12\, h_0(X)\,|\partial_t x(X,t)|^2\,dX\,dt.
\]
A standard variation $x\mapsto x+\eta$ with $\eta(\cdot,0)=\eta(\cdot,T)=0$
yields
\[
\delta \mathcal{A}
=
-\int_0^T\int_{\Omega_0} h_0(X)\,\partial_{tt}x(X,t)\,\eta(X,t)\,dX\,dt.
\]
Transforming back to Eulerian variables using $h_0(X)\,dX=h(x,t)\,dx$ and the
identity $\partial_{tt}x(X,t)=u_t(x,t)+u(x,t)u_x(x,t)$, we obtain the inertial
force 
\[
f_{\mathrm{inertial}}
=\frac{\delta \mathcal{A}}{\delta x}
=
-\,h\bigl(u_t+u u_x\bigr).
\]

We next derive the conservative force induced by the potential energy
\[
\mathcal{F}(h)
=
\int_\Omega \left(\frac12 g h^2 + g b(x)h\right)\,dx.
\]
Let $\xi(x,t)$ be the Eulerian displacement associated with the perturbed flow.
Mass conservation implies the transport identity
\begin{equation}
\label{eq:delta_h_transport}
\delta h = -\partial_x(h\,\xi).
\end{equation}
Using $\delta\mathcal F=\int (\delta\mathcal F/\delta h)\,\delta h\,dx$ with
\[
\frac{\delta \mathcal{F}}{\delta h}=g(h+b),
\]
we obtain
\[
\delta \mathcal{F}
=
-\int_\Omega g(h+b)\,\partial_x(h\,\xi)\,dx
=
\int_\Omega h\,\partial_x\!\bigl(g(h+b)\bigr)\,\xi\,dx,
\]
where the last step follows from integration by parts (assuming periodic
boundary conditions or vanishing boundary terms).
Hence,
\[
f_{\mathrm{conservative}}
=
\frac{\delta \mathcal{F}}{\delta x}
=
h\,\partial_x\!\bigl(g(h+b)\bigr).
\]
The EnVarA force balance
\[
\frac{\delta \mathcal{A}}{\delta x}
=
\frac{\delta \mathcal{F}}{\delta x},
\]
together with mass conservation \eqref{eq:SWE1}, recovers the momentum equation \eqref{eq:SWE2}.

A central advantage of EnVarA is that additional conservative effects can be
incorporated by augmenting the energy functional.
To introduce dispersion, we add the Fisher-information energy
\begin{equation}
\label{eq:fisher}
\mathcal{F}_{\mathrm{disp}}(h)
=
\frac{\varepsilon^2}{2}\int_\Omega \bigl|(\sqrt{h})_x\bigr|^2\,dx
=
\frac{\varepsilon^2}{8}\int_\Omega \frac{h_x^2}{h}\,dx,
\end{equation}
where $\varepsilon>0$ is the regularization parameter.
Energies of the form \eqref{eq:fisher} also appear in a different context in
Schr\"odinger bridge and entropic optimal transport problems, where the
Fisher-information term arises naturally from stochastic perturbations of the
continuity equation and can be interpreted as a noise-induced regularization;
see, e.g., \cite{Leonard2014,ChenGeorgiouPavon2016}.

A direct variational calculation yields the associated chemical potential
\begin{equation}
\label{eq:mu_from_fisher}
\frac{\delta \mathcal{F}_{\mathrm{disp}}}{\delta h}
=
-\frac{\varepsilon^2}{2}\,\frac{(\sqrt{h})_{xx}}{\sqrt{h}}
\;=:\;\mu,
\end{equation}
which is the \emph{quantum (Bohm) potential} in the Madelung literature
\cite{Madelung1927,CarlesDanchinSaut2012}.
Using the same transport identity \eqref{eq:delta_h_transport}, the dispersive
contribution to the conservative force is
\[
f_{\mathrm{disp}}
=
\frac{\delta \mathcal{F}_{\mathrm{disp}}}{\delta x}
=
h\,\mu_x.
\]
Combining the forces gives the regularized momentum balance
\begin{equation}
\label{disper}
h(u_t+u u_x)
+\,h\,\partial_x\!\bigl(g(h+b)\bigr) + h\,\mu_x=0.
\end{equation}
The regularized system \eqref{eq:SWE1} and \eqref{disper} is Hamiltonian and
satisfies the modified energy conservation law
\[
\frac{d}{dt}\left(\mathcal{E}(h,u) + \mathcal{F}_{\mathrm{disp}}(h)\right)=0,
\]
for smooth solutions and under appropriate boundary conditions.

% Alternatively, the Fisher-information regularization term \eqref{eq:fisher} can be naturally reformulated when we consider 
% the original energy law \eqref{ener0} with a modified kinetic relation as follows:
% \[
% h_t + (hu)_x - \varepsilon h_{xx} = 0.
% \]
% We leave out the detailed derivation in our future work. 

\subsection{Dispersive Regularization and the Nonlinear Schr\"odinger Equation}
The dispersively regularized shallow water system derived in the previous
subsection, cf.\ \eqref{eq:SWE1} and \eqref{disper}, admits an equivalent
formulation in terms of a nonlinear Schr\"odinger equation.
This equivalence follows from the classical Madelung transform and highlights
that the Fisher-information regularization naturally induces a
Schr\"odinger-type structure.

Specifically, assuming the existence of a real-valued velocity potential
$\phi$ such that $u=\phi_x$, and introducing the complex-valued wave function
\[
\psi(x,t)
=
\sqrt{h(x,t)}\,\exp\!\left(\frac{\mathrm{i}}{\varepsilon}\phi(x,t)\right),
\]
it can be shown that the regularized shallow water equations
\eqref{eq:SWE1}--\eqref{disper} are formally equivalent to the defocusing
nonlinear Schr\"odinger equation (NLS)
\begin{equation}
\label{eq:NLS}
\mathrm{i}\varepsilon \psi_t
=
-\frac{\varepsilon^2}{2}\psi_{xx}
+ g|\psi|^2\psi
+ g b(x)\psi,
\end{equation}
where the bottom topography enters as an external potential.

The hydrodynamic variables, namely the water height $h$ and the discharge
$q=hu$, are recovered from the wave function $\psi$ through the relations
\begin{equation}
\label{eq:recover_hq}
h = |\psi|^2,
\qquad
q
=
\varepsilon\,\mathrm{Im}\!\left(\overline{\psi}\,\psi_x\right).
\end{equation}

The main objective of this work is to numerically solve the nonlinear
Schr\"odinger equation \eqref{eq:NLS} in the small-$\varepsilon$ regime and to
use the resulting Madelung variables \eqref{eq:recover_hq} as an approximation
to the solution of the original dispersionless shallow water equations
\eqref{eq:SWE}.

\subsection{Semiclassical Limit and Approximation Regimes}

The convergence of solutions of the nonlinear Schr\"odinger equation
\eqref{eq:NLS} in the limit $\varepsilon \to 0$ has been extensively studied in
the context of the semiclassical (or zero-dispersion) limit of NLS; see, for
example, \cite{Grenier1998,Carles2008,CarlesDanchinSaut2012}.
In regimes where the limiting solution of the dispersionless shallow water
equations \eqref{eq:SWE1}--\eqref{eq:SWE2} remains smooth and the water height
$h$ stays uniformly bounded away from zero, it is well known that the Madelung
variables \eqref{eq:recover_hq} associated with the Schr\"odinger solution
converge locally in time to the corresponding strong solution of the shallow
water system.

When the dispersionless shallow water equations develop shocks (e.g., in
supercritical regimes), strong pointwise convergence is no longer expected.
In this case, solutions of the nonlinear Schr\"odinger equation typically
develop rapidly oscillatory wave trains, known as \emph{dispersive shock waves},
with characteristic wavelength of order $O(\varepsilon)$; see, for instance,
\cite{LaxLevermore1983,HoeferAblowitzEtAl2006,ElHoefer16}.
The limiting behavior in this regime is commonly described in a weak or averaged
sense using Whitham modulation theory \cite{Whitham74}.
As a consequence, the nonlinear Schr\"odinger equation is not expected to provide
a pointwise approximation of entropy solutions of the shallow water equations
after shock formation (often referred to as gradient catastrophe in the
semiclassical NLS literature).

The primary focus of the present work is the numerical reliability of the
Schr\"odinger-based approximation in shock-free regimes where wetting--drying
interfaces may occur.
The convergence theory in the presence of moving wetting--drying fronts remains
far less developed due to the loss of smoothness near vacuum boundaries.
Our numerical results indicate that, even in this setting, the nonlinear
Schr\"odinger equation provides accurate approximations of the shallow water
solution.
A salient feature of the proposed approach is that no ad hoc wetting--drying
treatment is required: the reconstructed water height $h=|\psi|^2$ is
guaranteed to remain non-negative throughout the computation.

\subsection{Numerical Scheme for NLS}
The numerical approximation of NLS \eqref{eq:NLS} has been extensively studied.
A wide variety of  discretization strategies have been proposed,
including finite difference methods, finite element/discontinuous
Galerkin methods, and spectral or pseudo-spectral methods; see, for example,
\cite{Akrivis92, BaoJinMarkowich2002, BaoJinMarkowich2002a, BaoCai2013,GuoXu2015}.

In the semiclassical regime $\varepsilon \ll 1$, the numerical approximation of
\eqref{eq:NLS} becomes particularly challenging due to the presence of
highly oscillatory solutions. It is well known that spatial mesh sizes and time step sizes must be chosen in
relation to the small parameter $\varepsilon$ in order to accurately resolve
the oscillatory behavior; see, e.g., \cite{MarkowichPietraPohl2002,BaoJinMarkowich2002a}.
Failure to respect these resolution constraints typically leads to loss of
accuracy or spurious numerical artifacts.

In this work, we adopt a spectral element method for the spatial discretization
of \eqref{eq:NLS}.
For the temporal discretization, we employ a second-order Strang splitting
scheme; see \cite{BaoJinMarkowich2002}, in which the linear dispersive operator is
separated from the nonlinear potential terms.
This strategy yields an efficient algorithm: the dispersive substep reduces to
the solution of a constant-coefficient linear problem, while the nonlinear
potential substep can be solved exactly and locally at each degree of
freedom.
As a result, the overall scheme is computationally efficient and particularly
well suited for simulations in the small-$\varepsilon$ regime.

For simplicity, we consider the one-dimensional spatial domain
$\Omega = [-L,L]$.
Let $\Omega_h=\{T\}$ be a conforming partition of $\Omega$ into elements.
We define the complex-valued spectral element space
\begin{equation}
\label{fes}
V_h^k
=
\bigl\{
v \in H^1(\Omega) :
v|_T \in \mathcal{P}_k(T), \ \forall T \in \Omega_h
\bigr\},
\end{equation}
where $k\ge1$ denotes the polynomial degree.
As is standard in spectral element methods, Gauss--Lobatto interpolation points
and basis functions are employed in the implementation.

We assume periodic boundary conditions for \eqref{eq:NLS}
(other boundary conditions will be discussed in the next subsection)
and define the periodic finite element space
\[
V_{h,\mathrm{per}}^k
:=
\bigl\{
v \in V_h^k : v(-L)=v(L)
\bigr\}.
\]

The spatial discretization of \eqref{eq:NLS} then reads:
find $\psi_h(t)\in V_{h,\mathrm{per}}^k$ such that
\begin{equation}
\label{spatial}
\bigl(
\mathrm{i}\varepsilon \psi_{h,t}, \phi_h
\bigr)_h
-
\frac{\varepsilon^2}{2}
\bigl(
\psi_{h,x}, \phi_{h,x}
\bigr)_h
=
\bigl(
g|\psi_h|^2\psi_h + g b(x)\psi_h,\phi_h
\bigr)_h,
\qquad
\forall \phi_h\in V_{h,\mathrm{per}}^k.
\end{equation}
Here $(\cdot,\cdot)_h$ denotes the discrete $L^2$ inner product obtained by
approximating integrals over each element using the $(k+1)$-point Gauss--Lobatto
quadrature rule.
This choice results in a diagonal mass matrix.

For the temporal discretization, we employ a second-order Strang splitting
method, in which \eqref{eq:NLS} is decomposed into linear dispersive and nonlinear
potential subproblems:
\begin{subequations}
\label{strang-split}
\begin{align}
\label{dis}
\mathrm{i}\varepsilon \psi_t &=
- \frac{\varepsilon^2}{2}\psi_{xx},
\\
\label{pot}
\mathrm{i}\varepsilon \psi_t &=
g|\psi|^2\psi + g b(x)\psi .
\end{align}
\end{subequations}

Starting from a known approximation $\psi_h^n\approx \psi(t_n)$, the Strang
splitting algorithm with time step size $\Delta t >0 $ proceeds as follows.

{\it Step 1: Half-step nonlinear potential update.}
We evolve \eqref{pot} over the half time interval $[t_n,t_n+\tfrac12\Delta t]$,
with initial condition $\psi_h(t_n)=\psi_h^n$.
The semi-discrete formulation reads
\begin{equation}
\label{potential}
\bigl(
\mathrm{i}\varepsilon \psi_{h,t}, \phi_h
\bigr)_h
=
\bigl(
g|\psi_h|^2\psi_h + g b(x)\psi_h,\phi_h
\bigr)_h,
\qquad
\forall \phi_h\in V_{h,\mathrm{per}}^k .
\end{equation}
Due to the use of Gauss--Lobatto basis functions and quadrature, this system
decouples at each degree of freedom, leading to the following nodal ODEs:
\begin{equation}
\label{potential-d}
\mathrm{i}\varepsilon \frac{\mathrm{d}}{\mathrm{d}t}\psi_h(x_j,t)
=
g\bigl(|\psi_h(x_j,t)|^2+b(x_j)\bigr)\psi_h(x_j,t),
\qquad
\forall x_j\in X_h,
\end{equation}
where $X_h=\{x_j\}$ denotes the set of all Gauss--Lobatto quadrature nodes on
$\Omega_h$.

Multiplying \eqref{potential-d} by $\overline{\psi_h}(x_j,t)$ and extracting the
imaginary part shows that
$\frac{\mathrm{d}}{\mathrm{d}t}|\psi_h(x_j,t)|^2=0$ for all $x_j \in X_h$.
Consequently, the nonlinear potential step preserves $|\psi_h|$ pointwise and
can be integrated exactly:
\begin{equation}
\label{step1}
\psi_h(x_j,t)
=
\exp\!\left(
-\frac{\mathrm{i}}{\varepsilon}
g\bigl(|\psi_h(x_j,t_n)|^2 + b(x_j)\bigr)(t-t_n)
\right)
\psi_h(x_j,t_n),
\qquad
x_j\in X_h.
\end{equation}
We denote the resulting solution at $t_n+\tfrac12\Delta t$ by
$\psi_h^{n,*}=\psi_h(t_n+\tfrac12\Delta t)$.

{\it Step 2: Full-step linear dispersive update.}
We next evolve the linear dispersive subproblem \eqref{dis} over the full time step
$\Delta t$, starting from $\psi_h^{n,*}$.
A Crank--Nicolson discretization in time yields
\begin{equation}
\label{dispersion}
\biggl(
\mathrm{i}\varepsilon
\frac{\psi_h^{n+1,*}-\psi_h^{n,*}}{\Delta t},
\phi_h
\biggr)_h
-
\frac{\varepsilon^2}{2}
\bigl(
(\psi_{h,x})^{n+1/2,*}, \phi_{h,x}
\bigr)_h
= 0,
\qquad
\forall \phi_h\in V_{h,\mathrm{per}}^k,
\end{equation}
where
\[
\psi_h^{n+1/2,*}
=
\frac12\bigl(\psi_h^{n+1,*}+\psi_h^{n,*}\bigr).
\]
This step requires the solution of a linear system with constant coefficients.
In our implementation, the resulting system is solved efficiently using a
preconditioned MINRES iterative solver equipped with a spectrally equivalent
block-diagonal preconditioner.

{\it Step 3: Half-step nonlinear potential update.}
Finally, we apply the nonlinear potential step \eqref{pot} again over the half
interval $[t_n+\tfrac12\Delta t,t_{n+1}]$, yielding
\begin{equation}
\label{step3}
\psi_h(x_j,t)
=
\exp\!\left(
-\frac{\mathrm{i}}{\varepsilon}
g\bigl(|\psi_h^{n+1,*}(x_j)|^2 + b(x_j)\bigr)
(t-t_n-\tfrac12\Delta t)
\right)
\psi_h^{n+1,*}(x_j),
\qquad
x_j\in X_h.
\end{equation}
The final approximation is then defined by
$\psi_h^{n+1}(x)=\psi_h(x,t_{n+1})\in V_{h,\mathrm{per}}^k$.

\subsection{Boundary Conditions}
\label{sec:bdry}

We conclude this section with a discussion of the treatment of far-field boundary
conditions for the nonlinear Schr\"odinger equation \eqref{eq:NLS}.
Many problems of interest are naturally posed on the whole space $\mathbb{R}$,
whereas numerical simulations must be carried out on a bounded computational
domain.
It is therefore necessary to impose suitable artificial boundary conditions on
a truncated domain that minimize spurious reflections and accurately represent
the far-field behavior of the solution.

In this work, we employ a sponge-layer (or absorbing-layer) approach to model
far-field boundary conditions.
The idea is to introduce a damping region near the artificial boundaries of the
computational domain, where outgoing waves are gradually attenuated before
reaching the boundary.
This strategy allows waves generated in the interior to exit the computational
domain with minimal reflection and provides a simple and robust treatment of
far-field boundaries for the Schr\"odinger equation.

Specifically, we extend the spatial domain from
$\Omega=[-L,L]$ to $\Omega_{\mathrm{ext}}=[-(L+\ell),\,L+\ell]$, where $\ell>0$
denotes the width of the sponge layer.
The regions $[-(L+\ell),-L]\cup[L,L+\ell]$ constitute the sponge layers.
We then augment \eqref{eq:NLS} with a complex absorbing potential (CAP) term,
leading to the modified equation
\begin{equation}
\label{absorb}
\mathrm{i}\varepsilon \psi_t
=
-\frac{\varepsilon^2}{2}\psi_{xx}
+ g|\psi|^2\psi
+ g b(x)\psi
- \mathrm{i}\sigma(x)\psi .
\end{equation}
Here $\sigma(x)\ge0$ is a spatially dependent damping function supported only
within the sponge layers.

In our implementation, $\sigma(x)$ is chosen using a quintic smoothstep ramp,
\begin{equation}
\label{sigma}
\sigma(x)
=
\sigma_{\max}
\left[
6\left(\frac{|x|-L}{\ell}\right)^5
-15\left(\frac{|x|-L}{\ell}\right)^4
+10\left(\frac{|x|-L}{\ell}\right)^3
\right],
\qquad
L<|x|<L+\ell,
\end{equation}
and $\sigma(x)=0$ for $|x|\le L$.

The sponge-layer thickness is taken as
\begin{equation}
\label{width}
\ell = N\,\frac{2\pi\varepsilon}{|\omega|},
\quad \text{with } N=16,
\end{equation}
where $\omega$ denotes the dominant outgoing wave number.
The maximum damping strength is chosen as
\begin{equation}
\label{sigmamax}
\sigma_{\max}
=
-\frac{2\,\varepsilon|\omega|}{\ell}\,\log(10^{-6}),
\end{equation}
corresponding to a target amplitude reduction of $10^{-6}$ across the sponge
layer.

With the sponge layer in place, we impose homogeneous Neumann boundary
conditions $\psi_x=0$ at the boundaries of the extended domain, i.e., at
$|x|=L+\ell$.
Since outgoing waves are strongly attenuated within the sponge region, this
simple boundary condition introduces negligible reflection into the interior
computational domain.

To incorporate the sponge layers into the Strang splitting scheme described in
the previous subsection, we proceed as follows.
The linear dispersive substep \eqref{dis} is unchanged.
On the interior region $[-L,L]$, where $\sigma(x)=0$, the nonlinear potential
updates remain identical to \eqref{step1}--\eqref{step3}.
The only modification occurs in the sponge layers, where the complex absorbing
potential term $-\mathrm{i}\sigma(x)\psi$ introduces an additional damping
effect. At the Gauss--Lobatto nodes $x_j\in X_h$, the half-step update \eqref{step1} in the
sponge layer takes the closed form
\begin{equation}
\label{cap_update}
\psi_h(x_j)
\ \mapsto\
\exp\!\left(
-\frac{\mathrm{i}}{\varepsilon}
g\bigl(|\psi_h(x_j)|^2+b(x_j)\bigr)\frac{\Delta t}{2}
\right)
\exp\!\left(
-\frac{\sigma(x_j)}{\varepsilon}\frac{\Delta t}{2}
\right)\psi_h(x_j),
\qquad x_j\in X_h.
\end{equation}
The same modification is applied in the final nonlinear half-step of the Strang
splitting scheme.

% \bibliographystyle{plain}
% \bibliography{references}

% \end{document}

% =========================
\section{Numerical Results}
% =========================
We present a series of numerical experiments to validate the proposed
Schr\"odinger-based approximation framework for the one-dimensional shallow
water equations.
All simulations are performed using the lowest-order
spectral element discretization with polynomial degree $k=1$.
A uniform spatial mesh with resolution $\Delta x = 0.05\,\varepsilon$ is used,
and the time step size is chosen as $\Delta t = \Delta x$.
The gravitational constant is fixed to $g=1$ throughout.
In the results presented below, we take $\varepsilon=0.01$, which is
representative of the small-dispersion regime considered in this work.
All numerical simulations are carried out using the MFEM finite element library \cite{mfem}.

\subsection{Riemann-Type Test Problems}

In this subsection, we consider several representative Riemann problems for the
shallow water equations \eqref{eq:SWE}.
In all cases, the initial data consist of two constant states separated at
$x=0$,
\begin{equation}
\label{eq:riemann_data}
(h(x,0),u(x,0))=
\begin{cases}
(h_L,u_L), & x<0,\\[0.6ex]
(h_R,u_R), & x>0,
\end{cases}
\end{equation}
with flat bottom topography $b(x)=0$.
Depending on the choice of left and right states, the corresponding exact
solutions may involve rarefaction waves, vacuum states, or shock formation
\cite{toro2024computational}.

To initialize the nonlinear Schr\"odinger equation \eqref{eq:NLS}, we employ a
uniform regularization strategy for all Riemann problems.
The initial water height is approximated by the smooth profile
\begin{equation}
\label{eq:h0_riemann}
h_0(x)
=
\frac{h_L+h_R}{2}+\frac{h_R-h_L}{2}\tanh(x/\delta),
\qquad
\delta=1.2\,\varepsilon,
\end{equation}
and the initial wave function is defined as
\begin{equation}
\label{eq:psi0_riemann}
\psi(x,0)=\sqrt{h_0(x)}\,\exp\!\left(\frac{\mathrm{i}}{\varepsilon}\phi_0(x)\right).
\end{equation}
The phase function $\phi_0(x)$ is chosen so that its derivative provides a smooth
approximation of the Riemann initial velocity,
\begin{equation}
\label{eq:phi0prime}
\phi_0'(x)
=
\frac{u_R+u_L}{2}
+
\frac{u_R-u_L}{2}\,
\tanh\!\left(\frac{x}{\delta}\right).
\end{equation}
Here we take 
\[ \phi_0(x) = \frac{u_R+u_L}{2}x +\frac{u_R-u_L}{2}\delta \Bigl( |x/\delta| +\log\bigl(1+\mathrm{e}^{-2|x|/\delta}\bigr) \Bigr). \]

Unless otherwise stated, computations are performed on $\Omega=[-2,2]$ with
homogeneous Neumann boundary conditions and final time $t=0.6$.
We now discuss three representative choices of Riemann data.

\subsubsection{Dam-Break Problem with a Dry Bed}
\label{dry_bed}
We first consider a dam-break problem with a dry bed to the right, corresponding
to
\[
(h_L,u_L)=(1.0,0),\qquad (h_R,u_R)=(0,0).
\]
The exact solution consists of a single rarefaction wave, with water height
given by \cite{toro2024computational}
\[
h(x,t)=
\begin{cases}
1, & x/t<-1,\\[0.3em]
\dfrac{(2-x/t)^2}{9}, & -1\le x/t\le 2,\\[0.3em]
0, & x/t>2.
\end{cases}
\]

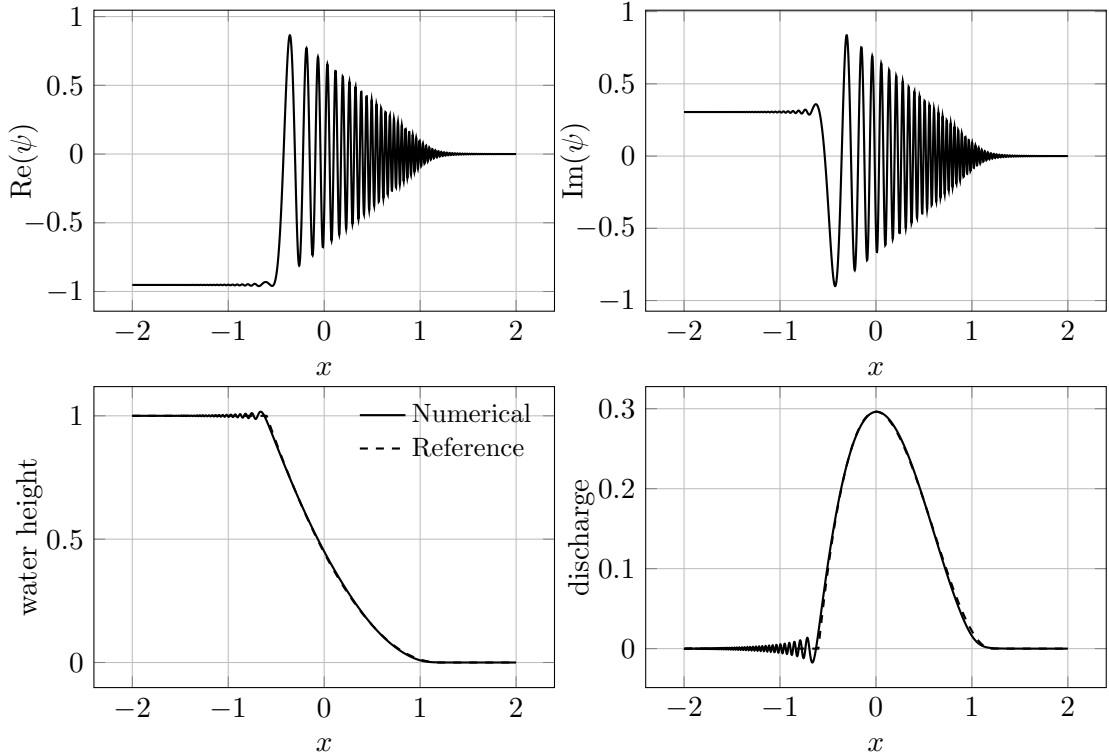
\begin{figure}[htbp]
\centering
\begin{tikzpicture}
\begin{groupplot}[
  group style={group size=2 by 2, horizontal sep=1.2cm, vertical sep=1.0cm},
  width=0.48\textwidth,
  height=0.35\textwidth,
  grid=both,
  xlabel={$x$},
  title = {Dam break (dry bed)},
  legend style={draw=none, fill=none, font=\small},
  legend cell align=left,
]

% ---------------- Panel (2,2) ----------------
\nextgroupplot[
  title={
  % Case 4 / time $t=t_4$
  },
        ylabel={$\mathrm{Re}(\psi)$},
    ylabel style={yshift=-13pt},
]
\addplot[thick] table[col sep=comma, x index=9, y index=5]{data/line_db0.csv};
% \addplot[thick, dashed] table[col sep=comma, x index=9, y index=7]{data/line_db0.csv};

% ---------------- Panel (2,1) ----------------
\nextgroupplot[
  title={
  % Case 3 / time $t=t_3$
  },
      ylabel={$\mathrm{Im}(\psi)$},
    ylabel style={yshift=-13pt},
]
\addplot[thick] table[col sep=comma, x index=9, y index=4]{data/line_db0.csv};
% \addplot[thick, dashed] table[col sep=comma, x index=9, y index=5]{data/line_db0.csv};

% ---------------- Panel (1,1) ----------------
\nextgroupplot[
  title={
  %Case 1 / time $t=t_1$
  },
  ylabel={water height},
  ylabel style={yshift=-6pt},
]
\addplot[thick] table[col sep=comma, x index=9, y index=0]{data/line_db0.csv};
\addlegendentry{Numerical}
\addplot[thick, dashed] table[col sep=comma, x index=9, y index=1]{data/line_db0.csv};
\addlegendentry{Reference}

% ---------------- Panel (1,2) ----------------
\nextgroupplot[
  title={
  %Case 2 / time $t=t_2$
  },
    ylabel={discharge},
    ylabel style={yshift=-6pt},
]
\addplot[thick] table[col sep=comma, x index=9, y index=2]{data/line_db0.csv};
\addplot[thick, dashed] table[col sep=comma, x index=9, y index=3]{data/line_db0.csv};

\end{groupplot}
\end{tikzpicture}
\caption{Numerical solution of the dam-break problem with a dry bed described in
 Subsection~\ref{dry_bed}.
The top row shows the real and imaginary parts of the Schr\"odinger wave function
$\psi$.
The bottom row displays the reconstructed hydrodynamic variables, namely the
water height $h=|\psi|^2$ and the discharge
$q=\varepsilon\,\mathrm{Im}(\overline{\psi}\,\psi_x)$.
The numerical solutions are compared with the corresponding
{dispersionless shallow water Riemann solution}.
}
\label{fig:ex1}
\end{figure}

Figure~\ref{fig:ex1} displays the computed Schr\"odinger wave function together with the
reconstructed water height and discharge at the final time $t=0.6$.
As expected in the small-$\varepsilon$ regime, the real and imaginary parts of
the wave function exhibit rapidly oscillatory behavior with characteristic
frequency of order $\mathcal{O}(1/\varepsilon)$.
Interestingly, these oscillations largely cancel when the solution is expressed
in terms of the reconstructed hydrodynamic variables.
In particular, both the water height and the discharge remain smooth across the
wet--dry interface and agree very well with the corresponding
\emph{dispersionless shallow water Riemann solution} in this region.

Small-amplitude oscillations of order $\mathcal{O}(\varepsilon)$ are also
observed near the head of the left-going rarefaction wave.
Such oscillations are consistent with dispersive effects introduced by the
Schr\"odinger regularization and diminish as $\varepsilon$ decreases.
Overall, the results demonstrate that, despite the highly oscillatory nature of
the underlying wave function, the Schr\"odinger-based approach provides an
accurate approximation of the shallow water solution in the shock-free regime,
including in the presence of a moving wet--dry interface.

\subsubsection{Dam-Break Problem with a Wet Bed}
\label{wet_bed}

We next consider a dam-break Riemann problem with a wet bed on the right, which
leads to shock formation.
The initial left and right states are given by
\[
(h_L,u_L)=(1.0,0),\qquad (h_R,u_R)=(0.2,0).
\]
For this configuration, the exact solution of the dispersionless shallow water
equations consists of a left-going rarefaction wave followed by a right-going
shock \cite{toro2024computational}.

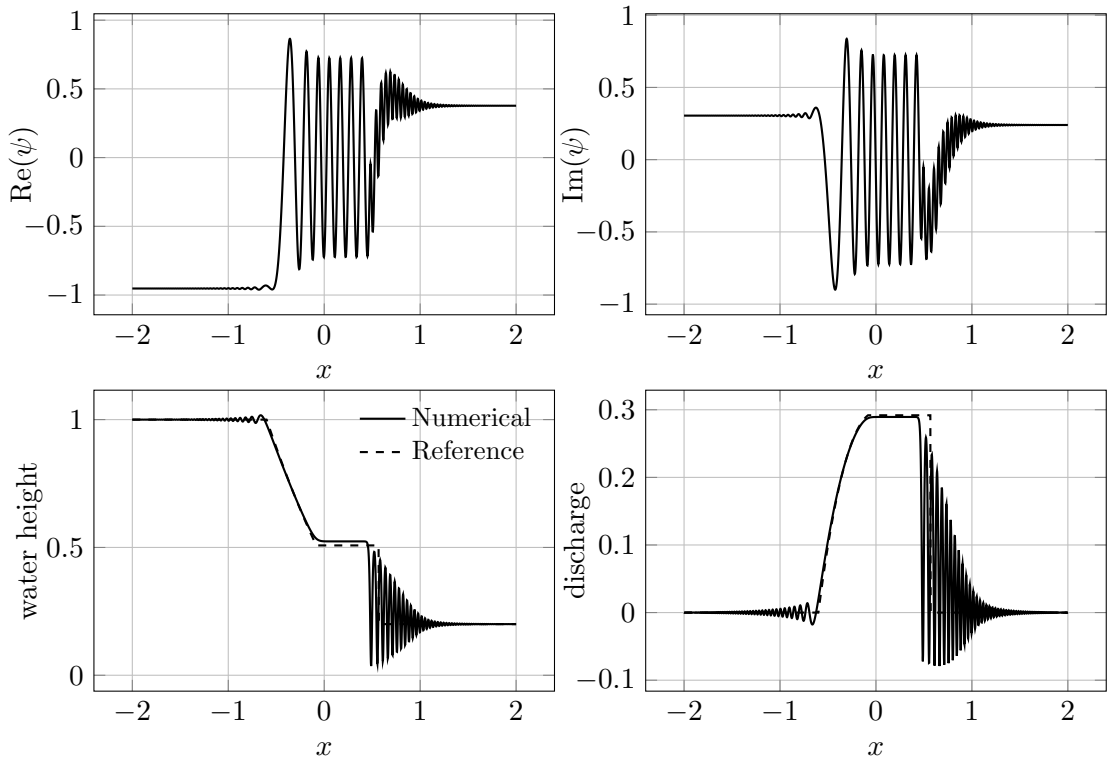
\begin{figure}[htbp]
\centering
\begin{tikzpicture}
\begin{groupplot}[
  group style={group size=2 by 2, horizontal sep=1.2cm, vertical sep=1.0cm},
  width=0.48\textwidth,
  height=0.35\textwidth,
  grid=both,
  xlabel={$x$},
  title = {Dam break (dry bed)},
  legend style={draw=none, fill=none, font=\small},
  legend cell align=left,
]

% ---------------- Panel (2,2) ----------------
\nextgroupplot[
  title={
  % Case 4 / time $t=t_4$
  },
        ylabel={$\mathrm{Re}(\psi)$},
    ylabel style={yshift=-13pt},
]
\addplot[thick] table[col sep=comma, x index=9, y index=5]{data/line_db1.csv};

% ---------------- Panel (2,1) ----------------
\nextgroupplot[
  title={
  % Case 3 / time $t=t_3$
  },
      ylabel={$\mathrm{Im}(\psi)$},
    ylabel style={yshift=-13pt},
]
\addplot[thick] table[col sep=comma, x index=9, y index=4]{data/line_db1.csv};

% ---------------- Panel (1,1) ----------------
\nextgroupplot[
  title={
  %Case 1 / time $t=t_1$
  },
  ylabel={water height},
  ylabel style={yshift=-6pt},
]
\addplot[thick] table[col sep=comma, x index=9, y index=0]{data/line_db1.csv};
\addlegendentry{Numerical}
\addplot[thick, dashed] table[col sep=comma, x index=9, y index=1]{data/line_db1.csv};
\addlegendentry{Reference}

% ---------------- Panel (1,2) ----------------
\nextgroupplot[
  title={
  %Case 2 / time $t=t_2$
  },
    ylabel={discharge},
    ylabel style={yshift=-13pt},
]
\addplot[thick] table[col sep=comma, x index=9, y index=2]{data/line_db1.csv};
\addplot[thick, dashed] table[col sep=comma, x index=9, y index=3]{data/line_db1.csv};

\end{groupplot}
\end{tikzpicture}
\caption{Numerical solution of the dam-break problem with a wet bed described in
 Subsection~\ref{wet_bed}.
The top row shows the real and imaginary parts of the Schr\"odinger wave function
$\psi$.
The bottom row displays the reconstructed hydrodynamic variables, namely the
water height $h=|\psi|^2$ and the discharge
$q=\varepsilon\,\mathrm{Im}(\overline{\psi}\,\psi_x)$.
The numerical solutions are compared with the corresponding
{dispersionless shallow water Riemann solution}.
}
\label{fig:ex2}
\end{figure}

Figure~\ref{fig:ex2} shows the numerical solution of the dam-break problem with a wet bed at
the final time $t=0.6$.
As in the previous example, the top row displays the real and imaginary parts of
the Schr\"odinger wave function, which exhibit rapidly oscillatory behavior in
the small-$\varepsilon$ regime.
The bottom row presents the reconstructed hydrodynamic variables, namely the
water height and the discharge, together with the corresponding dispersionless
shallow water Riemann solution.

In the region associated with the left-going rarefaction wave, the reconstructed
water height and discharge agree well with the dispersionless shallow water
solution. In this regime, the oscillations present in the wave function largely
cancel when expressed in terms of the hydrodynamic variables, and the
Schr\"odinger-based approximation provides a reasonable representation of the
underlying shallow water dynamics.

In contrast, near the right-going shock, the Schr\"odinger solution no longer
approximates the dispersionless shock profile.
Instead of converging to a monotone entropy shock as in the shallow water
equations, the solution develops a rapidly oscillatory structure characteristic
of a dispersive shock wave.
These oscillations persist in the reconstructed water height and discharge and
do not vanish as $\varepsilon$ decreases.

This behavior reflects a fundamental limitation of the present approach: the
nonlinear Schr\"odinger regularization replaces gradient catastrophe in the
dispersionless shallow water equations by dispersive shock waves rather than
classical entropy shocks.
Consequently, while the Schr\"odinger-based formulation performs well in
shock-free regimes, it is not suitable for approximating shallow water solutions
in the presence of shock formation.

\subsubsection{Generation of Vacuum from Wet Initial States}
\label{vac}
Finally, we consider a Riemann problem in which vacuum is generated dynamically.
The initial states are
\[
(h_L,u_L)=(1.0,-3.0),\qquad (h_R,u_R)=(2.0,3.0).
\]
Since
\[
2(a_L+a_R)<u_R-u_L,
\qquad
a_L=\sqrt{g h_L},\quad a_R=\sqrt{g h_R},
\]
two rarefaction waves propagate outward, creating an expanding dry region
between them.
The exact solution contains a vacuum state and is known analytically
\cite{toro2024computational}.

Because the initial velocity is nonzero, we employ the sponge-layer boundary
treatment described in Subsection~\ref{sec:bdry}.
The computational domain is extended to
$\Omega_{\mathrm{ext}}=[-(L+\ell),L+\ell]$ with $L=2$ and
$\ell=16(2\pi\varepsilon/3)\approx 0.335$.

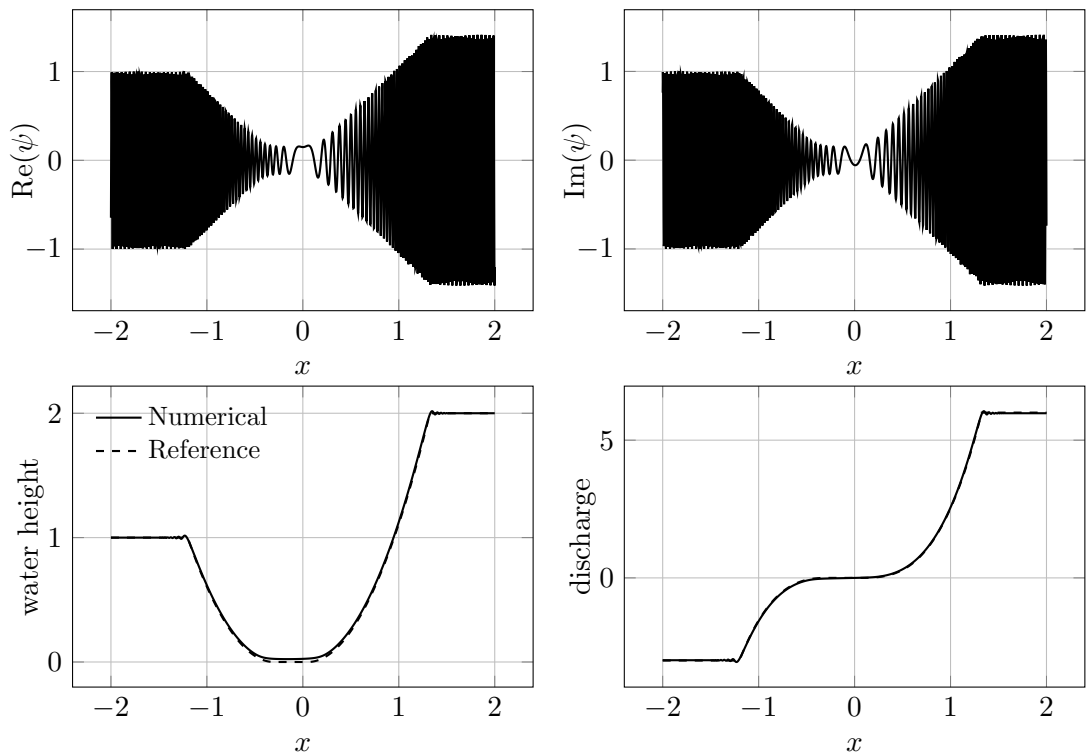
\begin{figure}[htbp]
\centering
\begin{tikzpicture}
\begin{groupplot}[
  group style={group size=2 by 2, horizontal sep=1.2cm, vertical sep=1.0cm},
  width=0.48\textwidth,
  height=0.35\textwidth,
  grid=both,
  xlabel={$x$},
  title = {Dam break (dry bed)},
  legend style={draw=none, fill=none, font=\small},
  legend cell align=left,
]

% ---------------- Panel (2,2) ----------------
\nextgroupplot[
  title={
  % Case 4 / time $t=t_4$
  },
        ylabel={$\mathrm{Re}(\psi)$},
    ylabel style={yshift=-13pt},
]
\addplot[thick] table[col sep=comma, x index=9, y index=5]{data/line_drying0.csv};

% ---------------- Panel (2,1) ----------------
\nextgroupplot[
  title={
  % Case 3 / time $t=t_3$
  },
      ylabel={$\mathrm{Im}(\psi)$},
    ylabel style={yshift=-13pt},
]
\addplot[thick] table[col sep=comma, x index=9, y index=4]{data/line_drying0.csv};

% ---------------- Panel (1,1) ----------------
\nextgroupplot[
  title={
  %Case 1 / time $t=t_1$
  },
  ylabel={water height},
  ylabel style={yshift=-6pt},
    legend pos=north west,
]
\addplot[thick] table[col sep=comma, x index=9, y index=0]{data/line_drying0.csv};
\addlegendentry{Numerical}
\addplot[thick, dashed] table[col sep=comma, x index=9, y index=1]{data/line_drying0.csv};
\addlegendentry{Reference}

% ---------------- Panel (1,2) ----------------
\nextgroupplot[
  title={
  %Case 2 / time $t=t_2$
  },
    ylabel={discharge},
    ylabel style={yshift=-6pt},
]
\addplot[thick] table[col sep=comma, x index=9, y index=2]{data/line_drying0.csv};
\addplot[thick, dashed] table[col sep=comma, x index=9, y index=3]{data/line_drying0.csv};

\end{groupplot}

\end{tikzpicture}
\caption{Numerical solution of the vacuum generation problem described in
 Subsection~\ref{vac}.
The top row shows the real and imaginary parts of the Schr\"odinger wave function
$\psi$.
The bottom row displays the reconstructed hydrodynamic variables, namely the
water height $h=|\psi|^2$ and the discharge
$q=\varepsilon\,\mathrm{Im}(\overline{\psi}\,\psi_x)$.
The numerical solutions are compared with the corresponding
{dispersionless shallow water Riemann solution}.
}
\label{fig:ex3}
\end{figure}

Figure~\ref{fig:ex3} shows the numerical solution of the vacuum generation problem
at the final time $t=0.3$.
As in the previous examples, the Schr\"odinger wave function exhibits highly
oscillatory behavior.
In particular, strong oscillations are observed near the left and right
far-field regions, where the amplitude approaches constant states and the phase
of the wave function varies rapidly.
At the scale of the full domain, these oscillations are too rapid to be resolved
visually and therefore appear saturated in the plots of
$\mathrm{Re}(\psi)$ and $\mathrm{Im}(\psi)$.

Despite this highly oscillatory behavior at the level of the wave function, the
reconstructed hydrodynamic variables remain well behaved.
Both the water height and the discharge are smooth across the rarefaction waves
and accurately capture the formation of the central dry region.
In particular, the numerical solution agrees very well with the dispersionless
shallow water Riemann solution.

Small-amplitude oscillations of order $\mathcal{O}(\varepsilon)$ are visible in
the reconstructed variables near the edges of the rarefaction fans.
These oscillations are consistent with dispersive effects introduced by the
Schr\"odinger regularization and diminish as $\varepsilon$ decreases.
Overall, this example further demonstrates that, although the Schr\"odinger wave
function itself is highly oscillatory, the Schr\"odinger-based formulation
provides an accurate approximation of the shallow water dynamics in regimes
where the solution remains free of shock formation, including cases where a
vacuum state emerges dynamically.

\subsection{Oscillating lake problem over a parabolic bowl}
\label{bowl}
We next consider the classical oscillating lake problem over a parabolic bowl
topography \cite{Thacker81,Sampson05}, formulated here in a nondimensionalized
setting.
The bottom topography is given by
\[
b(x)=x^2,
\]
with initial free-surface elevation, water height, and velocity
\[
\eta_0(x)=\max\{0.5-\sqrt{2}\,x, b(x)\},
\qquad
h_0(x)=\eta_0(x)-b(x),
\qquad
u_0(x)=0.
\]
This problem admits an exact periodic solution with moving wetting--drying
interfaces,
\[
h(x,t)=\eta(x,t)-b(x),
\]
where
\[
\eta(x,t)=\max\{0.75-0.25\cos(2\sqrt{2}\,t)-\sqrt{2}\,x\cos(\sqrt{2}\,t), b(x)\}.
\]

The computational domain is taken as $\Omega=[-2,2]$, equipped with homogeneous
Neumann boundary conditions.
The initial Schr\"odinger wave function is constructed from a regularized free
surface profile,
\[
\psi(x,0)=
\sqrt{\delta\log\!\left(1+\exp\!\left(\frac{\eta_0(x)-b(x)}{\delta}\right)\right)},
\qquad
\delta=1.2\,\varepsilon.
\]

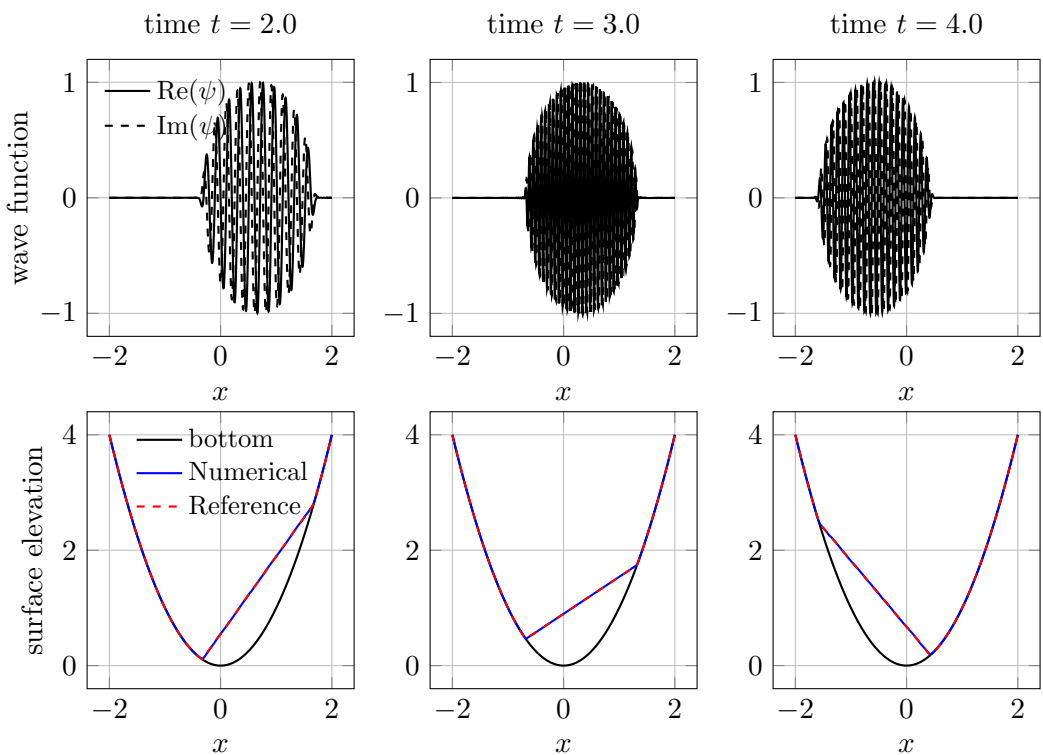
\begin{figure}[htbp]
\centering
\begin{tikzpicture}
\begin{groupplot}[
  group style={group size=3 by 2, horizontal sep=1.0cm, vertical sep=1.0cm},
  width=0.32\textwidth,
  height=0.33\textwidth,
  grid=both,
  xlabel={$x$},
  legend style={draw=none, fill=none, font=\small},
  legend cell align=left,
]

% ======================================================
% Top row: wave function psi
% ======================================================

% --- t = 2.0
\nextgroupplot[
  title={time $t=2.0$},
  ylabel={wave function},
  ylabel style={yshift=-3pt},
  legend pos=north west,
]
\addplot[thick] table[col sep=comma, x index=10, y index=5]{data/line_lake2.csv};
\addlegendentry{$\mathrm{Re}(\psi)$}
\addplot[thick, dashed] table[col sep=comma, x index=10, y index=6]{data/line_lake2.csv};
\addlegendentry{$\mathrm{Im}(\psi)$}

% --- t = 3.0 (NEW)
\nextgroupplot[
  title={time $t=3.0$},
]
\addplot[thick] table[col sep=comma, x index=10, y index=5]{data/line_lake3.csv};
\addplot[thick, dashed] table[col sep=comma, x index=10, y index=6]{data/line_lake3.csv};

% --- t = 4.0
\nextgroupplot[
  title={time $t=4.0$},
]
\addplot[thick] table[col sep=comma, x index=10, y index=5]{data/line_lake4.csv};
\addplot[thick, dashed] table[col sep=comma, x index=10, y index=6]{data/line_lake4.csv};

% ======================================================
% Bottom row: surface elevation eta
% ======================================================

% --- t = 2.0
\nextgroupplot[
  ylabel={surface elevation},
]
\addplot[thick] table[col sep=comma, x index=10, y index=0]{data/line_lake2.csv};
\addplot[thick, blue] table[col sep=comma, x index=10, y index=2]{data/line_lake2.csv};
\addplot[thick, dashed, red] table[col sep=comma, x index=10, y index=1]{data/line_lake2.csv};
\addlegendentry{bottom}
\addlegendentry{Numerical\quad\quad}
\addlegendentry{Reference}

% --- t = 3.0 (NEW)
\nextgroupplot[
]
\addplot[thick] table[col sep=comma, x index=10, y index=0]{data/line_lake3.csv};
\addplot[thick, blue] table[col sep=comma, x index=10, y index=2]{data/line_lake3.csv};
\addplot[thick, dashed, red] table[col sep=comma, x index=10, y index=1]{data/line_lake3.csv};

% --- t = 4.0
\nextgroupplot[
]
\addplot[thick] table[col sep=comma, x index=10, y index=0]{data/line_lake4.csv};
\addplot[thick, blue] table[col sep=comma, x index=10, y index=2]{data/line_lake4.csv};
\addplot[thick, dashed, red] table[col sep=comma, x index=10, y index=1]{data/line_lake4.csv};

\end{groupplot}
\end{tikzpicture}

\caption{Numerical solution of the oscillating lake problem described in
Subsection~\ref{bowl}.
The top row shows the real and imaginary parts of the Schr\"odinger wave function
$\psi$.
The bottom row displays the reconstructed free-surface elevation
$\eta(x)=h(x)+b(x)$, together with the bottom topography $b(x)$.
The left, middle, and right columns correspond to times $t=2.0$, $t=3.0$, and
$t=4.0$, respectively.}
\label{fig:ex4}
\end{figure}

Figure~\ref{fig:ex4} presents the numerical solution at times $t=2.0$, $t=3.0$,
and $t=4.0$.
As in the previous examples, the real and imaginary parts of the Schr\"odinger
wave function exhibit rapidly oscillatory behavior in the small-$\varepsilon$
regime.
We observe that the oscillation frequency varies with time, reflecting the
time-dependent phase structure of the underlying wave function.

Despite these pronounced oscillations at the level of the wave function, the
reconstructed free-surface elevation remains smooth and closely follows the
exact solution of the dispersionless shallow water equations.
In particular, the numerical solution accurately captures the periodic motion
of the free surface as well as the evolution of the wetting--drying interfaces.
This example further indicates that the Schr\"odinger-based formulation is
applicable to this setting and is able to capture 
the evolution of wetting--drying interfaces without requiring any ad hoc
numerical treatment near dry regions.

\subsection{Well-Balanced Tests}
\label{wb}
Finally, we examine the well-balanced behavior of the Schr\"odinger-based
approximation for steady lake-at-rest configurations.
We consider the bottom topography
\[
b(x)=b_{\max}\mathrm{e}^{-10x^2},
\]
with initial water height and velocity
\[
h_0(x)=(1-b(x))_+,
\qquad
u_0(x)=0.
\]
The corresponding exact solution of the dispersionless shallow water equations
is the steady lake-at-rest state
\[
h(x,t)=(1-b(x))_+,
\qquad
u(x,t)=0.
\]

\begin{figure}[htbp]
\centering
\begin{tikzpicture}
\begin{groupplot}[
  group style={group size=2 by 3, horizontal sep=1.2cm, vertical sep=1.0cm},
  width=0.48\textwidth,
  height=0.35\textwidth,
  grid=both,
  xlabel={$x$},
  title = {Dam break (dry bed)},
  legend style={draw=none, fill=none, font=\small},
  legend cell align=left,
]

\nextgroupplot[
  title={
 $b_{\max}=0.9$
  },
      ylabel={wave function},
    ylabel style={yshift=0pt},
   ymin=-0.03,
  ymax=0.92,
 legend style={
    at={(0.7,0.4)},
    anchor=north west
  },
]
\addplot[thick] table[col sep=comma, x index=10, y index=5]{data/line_wb0.csv};
\addlegendentry{$\mathrm{Re}(\psi)$}
\addplot[thick, dashed] table[col sep=comma, x index=10, y index=6]{data/line_wb0.csv};
\addlegendentry{$\mathrm{Im}(\psi)$}

\nextgroupplot[
  title={
 $b_{\max}=1.1$
  },
     ymin=-0.03,
  ymax=0.92,
]
\addplot[thick] table[col sep=comma, x index=10, y index=5]{data/line_wb1.csv};
\addplot[thick, dashed] table[col sep=comma, x index=10, y index=6]{data/line_wb1.csv};

% ---------------- Panel (2,2) ----------------
\nextgroupplot[
  title={
  },
        ylabel={surface elevation},
    ylabel style={yshift=0pt},
      ymin=-0.1,
  ymax=1.2,
 legend style={
    at={(0.6,0.5)},
    anchor=north west
  }            
]
\addplot[thick] table[col sep=comma, x index=10, y index=0]{data/line_wb0.csv};
\addplot[thick, blue] table[col sep=comma, x index=10, y index=2]{data/line_wb0.csv};
\addplot[thick, dashed, red] table[col sep=comma, x index=10, y index=1]{data/line_wb0.csv};
\addlegendentry{bottom}
\addlegendentry{Numerical\qquad\qquad\qquad}
\addlegendentry{Reference}

% ---------------- Panel (2,1) ----------------
\nextgroupplot[
  title={
  % Case 4 / time $t=t_4$
  },
    ymin=-0.1,
  ymax=1.2,
]
\addplot[thick] table[col sep=comma, x index=10, y index=0]{data/line_wb1.csv};
\addplot[thick, blue] table[col sep=comma, x index=10, y index=2]{data/line_wb1.csv};
\addplot[thick, dashed, red] table[col sep=comma, x index=10, y index=1]{data/line_wb1.csv};

% ---------------- Panel (2,2) ----------------
\nextgroupplot[
  title={
  },
        ylabel={surface elevation},
    ylabel style={yshift=0pt},
  ymin=0.989,
  ymax=1.011,
  ytick={0.99,1.00,1.01},
  yticklabel style={
    /pgf/number format/fixed,
    /pgf/number format/precision=2
  },
 legend style={
    at={(0.6,0.8)},
    anchor=north west
  }            
]
\addplot[thick, blue] table[col sep=comma, x index=10, y index=2]{data/line_wb0.csv};
\addplot[thick, dashed, red] table[col sep=comma, x index=10, y index=1]{data/line_wb0.csv};
\addlegendentry{Numerical\qquad\qquad\qquad}
\addlegendentry{Reference}

% ---------------- Panel (2,1) ----------------
\nextgroupplot[
  title={
  % Case 4 / time $t=t_4$
  },
]
\addplot[thick, blue] table[col sep=comma, x index=10, y index=2]{data/line_wb1.csv};
\addplot[thick, dashed, red] table[col sep=comma, x index=10, y index=1]{data/line_wb1.csv};

\end{groupplot}

\end{tikzpicture}
\caption{Numerical solution at time $t= 1.0$ of the well-balanced problems described in
 Subsection~\ref{wb}.
The top row shows the real and imaginary parts of the Schr\"odinger wave function
$\psi$.
The middle row displays the reconstructed free-surface elevation
$\eta(x)=h(x)+b(x)$, where $h=|\psi|^2$, together with the bottom topography
$b(x)$.
The bottom row shows the reconstructed free-surface elevation
$\eta(x)=h(x)+b(x)$ without the bottom topography, in order to highlight
the small deviations of the numerical solution from the exact (dispersionless) steady state.
The left column corresponds to the fully wet case with $b_{\max}=0.9$, while the
right column corresponds to the partially dry case with $b_{\max}=1.1$.
}
\label{fig:ex5}
\end{figure}
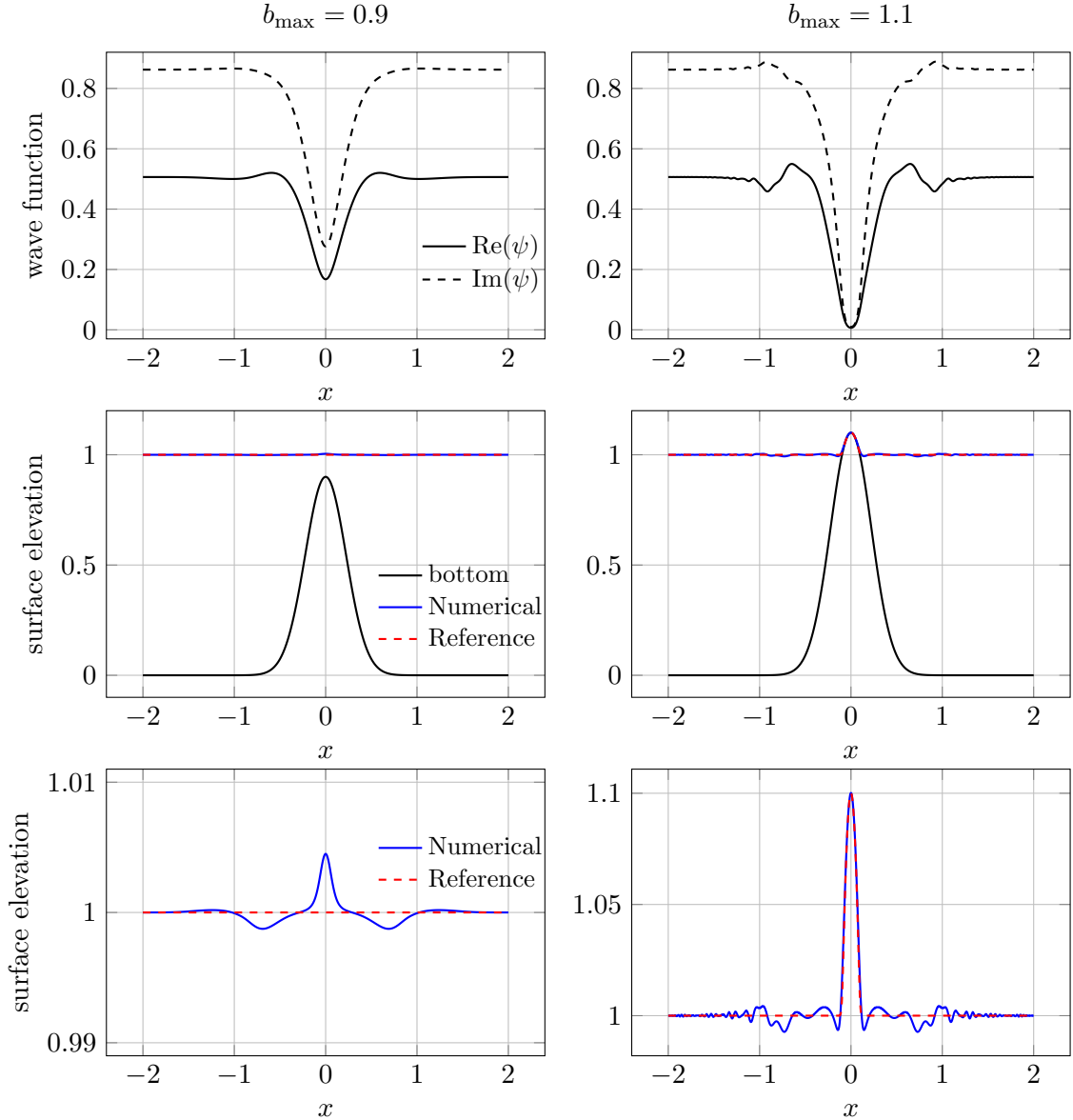

We impose periodic boundary conditions on $\Omega=[-2,2]$ and consider two
representative cases: $b_{\max}=0.9$, which corresponds to a fully wet steady
state, and $b_{\max}=1.1$, which produces a partially dry configuration.
Numerical results at the final time $t=1.0$ are shown in
Figure~\ref{fig:ex5}.

In contrast to the previous test cases, the Schr\"odinger wave
function exhibits significantly reduced oscillatory behavior in this setting.
This is consistent with the fact that the underlying shallow water solution is a
steady state with zero velocity, so that no rapidly varying phase is generated
in the Schr\"odinger formulation.
As a result, the real and imaginary parts of $\psi$ remain relatively smooth in
both the fully wet and partially dry cases.

Small deviations from the exact steady free-surface elevation are nevertheless
observed in the numerical solution.
In both cases, these deviations are of magnitude $\mathcal{O}(\varepsilon)$,
reflecting the fact that the nonlinear Schr\"odinger equation
\eqref{eq:NLS} does not satisfy an exact well-balanced property.
We also observe that the partially dry configuration exhibits slightly larger
approximation errors than the fully wet case, although the overall error level
remains of order $\mathcal{O}(\varepsilon)$ in both settings.

% =========================
\section{Conclusion}
% =========================

We have introduced a Schr\"odinger-based dispersive regularization framework for
the numerical simulation of the one-dimensional shallow water equations.
By solving a semilinear complex-valued Schr\"odinger equation and recovering the
hydrodynamic variables through postprocessing, the proposed approach transforms
the original nonlinear hyperbolic system into a form that can be efficiently
approximated using time-splitting and high-order spatial discretizations.
A key advantage of this formulation is its ability to naturally accommodate
vacuum states and moving wetting--drying interfaces without requiring ad hoc
numerical treatments.

Numerical experiments demonstrate that, in subcritical and shock-free regimes,
the Schr\"odinger-based approximation provides an accurate representation of the
dispersionless shallow water solution, with errors typically of order
$\mathcal{O}(\varepsilon)$.
This includes challenging scenarios involving wetting and drying, vacuum
generation, and time-dependent free-surface motion over nontrivial bathymetry.
In contrast, when shock formation occurs in the shallow water equations, the
Schr\"odinger regularization produces dispersive shock waves rather than
classical entropy shocks, highlighting an inherent limitation of the present
approach in such regimes.

Overall, the results indicate that Schr\"odinger-based dispersive
regularization offers a promising alternative framework for the numerical
simulation of shallow water flows in regimes where solutions remain smooth or
feature vacuum and moving shorelines.
Future work will focus on rigorous error analysis in the presence of vacuum,
extensions to higher-dimensional settings, and applications to more realistic
geophysical and engineering flow problems.

\bibliographystyle{plain}
\bibliography{references}

\end{document}